\documentclass[12pt]{article}
\usepackage{amsthm,amstext,amsmath,amscd,amssymb,eucal,latexsym}
\usepackage[matrix,arrow]{xy} 

\title
{Construction of Rational Surfaces of Degree 12 \\
 in Projective Fourspace }
\author{Hirotachi Abo \\
}

\date{}

\newcommand{\C}{{\mathbb C}}

\newcommand{\F}{{\mathbb F}}
\newcommand{\GG}{{\mathbb G}}
\newcommand{\M}{{\mathbb M}}

\renewcommand{\P}{{\mathbb P}}
\newcommand{\Q}{{\mathbb Q}}

\newcommand{\Z}{{\mathbb Z}}

\newcommand{\LL}{{\mathbb L}}

\newcommand{\PP}{{\mathbb P}}

\newcommand{\Syz}{{\rm{Syz}\;}}

\newcommand{\h}{{\rm{h}}}
\newcommand{\HH}{{\rm{H}}}

\newcommand{\cO}{{\mathcal O}}

\newcommand{\s}{\mathcal}

\newcommand{\sE}{{\s E}}
\newcommand{\sF}{{\s F}}

\newcommand{\sI}{{\s I}}

\newcommand{\sX}{{\s X}}

\newcommand{\lra}{\longrightarrow}

\newcommand{\tensor}{\otimes}

\newcommand{\punkt}{\HHspace{-.3ex}\raise.15ex\HHbox to1ex{\HHuge.}}

\newlength{\br}
\newlength{\ho}

\DeclareMathOperator{\Spec}{Spec}
\DeclareMathOperator{\HHom}{Hom}

\DeclareMathOperator{\Ker}{Ker}
\DeclareMathOperator{\im}{im}

\DeclareMathOperator{\Coker}{Coker}

\DeclareMathOperator{\codim}{codim}

\newcommand{\integer}{\Z}

\newcommand{\Mac}{{\texttt {MACAULAY2}}}

\newcommand{\paper}{: \begin{it}}
\newcommand{\jour }{, \end{it}}

\newtheorem{theorem}{Theorem}[section]
\newtheorem{lemma}[theorem]{Lemma}
\newtheorem{proposition}[theorem]{Proposition}

\theoremstyle{definition}

\newtheorem{remark}[theorem]{Remark}

\numberwithin{equation}{section}

\begin{document}

\maketitle
\begin{abstract}
The aim of this paper is to present two different constructions 
of smooth rational surfaces in projective fourspace
with degree $12$ and sectional genus $13$. 
In particular, we establish the existences of five different families of 
smooth rational surfaces in projective fourspace with the prescribed 
invariants. 
\end{abstract}
\section{Introduction}
Hartshone conjectured that only finitely many components of the Hilbert 
scheme of surfaces in $\P^4$ contain smooth rational surfaces. 
In 1989, this conjecture was positively solved by Ellingsrud and Peskine~\cite{ep}. 
The exact bound for the degree is, however, still open, and  
hence the question concerning the exact bound motivates us to find smooth rational surfaces  
in $\P^4$.  The goal of this paper is to construct five different families 
of smooth rational surfaces in $\P^4$ with degree $12$ and sectional genus $13$. 
The rational surfaces in $\P^4$ were previously known up to degree $11$. 

In this paper, we present two different constructions of smooth rational surfaces in $\P^4$ 
with the given invariants. 
Both the constructions are based on the technique of ``Beilinson monad". 
Let $V$ be a finite-dimensional vector space over a field $K$ 
and let $W$ be its dual space.  
The basic idea behind a Beilinson monad is to represent 
a given coherent sheaf on $\P(W)$ as a homology of a finite complex 
whose objects are direct sums of  bundles of differentials. 
The differentials in the monad are given by homogeneous 
matrices over an exterior algebra $E=\bigwedge V$. 
To construct a Beilinson monad for a given coherent sheaf, 
we typically take the following steps:
Determine the type of the Beilinson monad, that is, 
determine each object, and then 
find differentials in the monad. 

Let $X$ be a smooth rational surface in $\P^4=\P(W)$ with degree $12$ 
and sectional genus $13$.  
The type of a Beilinson monad for the (suitably twisted) ideal sheaf of $X$ can be 
derived from the knowledge of its cohomology groups. 
Such information is partially determined from general 
results such as the Riemann-Roch formula and Kodaira vanishing theorem. 
It is, however, hard to determine the dimensions of all cohomology 
groups needed to determine the type 
of the Beilinson monad. For this reason,  
we assume that the ideal sheaf 
of $X$ has the so-called ``natural cohomology" in some range of twists.  
From this assumption, the Beilinson monad for the twisted 
ideal sheaf $\sI_X(4)$ of $X$ has the following form:
\begin{eqnarray}
 4\Omega^3(3) \stackrel{A}{\rightarrow} 2\Omega^2(2) \oplus 2\Omega^1(1) 
 \stackrel{B}{\rightarrow} 3\cO. 
 \label{eq:minimal-monad}
\end{eqnarray}
To detect differentials in (\ref{eq:minimal-monad}), we use the following techniques developed recently: (1) the fisrt technique is an exterior algebra method  
due to Eisenbud, Fl{\o}ystad and Schreyer~\cite{efs} 
and (2) the other one is the method using small finite fields and random trials due to Schreyer~\cite{schreyer}.  

\vspace{1mm}
(1) Eisenbud, Fl{\o}ystad and Schreyer  presented 
an explicit version of  the Bernstein-Gel'fand-Gel'fand correspondence. 
This correspondence is an isomorphism 
between the derived category of bounded complexes of  
finitely generated $S$-graded modules and the derived category 
of certain ``Tate resolutions" of $E$-modules, 
where $S=\mathrm{Sym}_K(W)$. 
As an application, they constructed the Beilinson monad from 
the Tate resolution explicitly. 
This enables us to describe the conditions that 
the differentials in the Beilinson monad must satisfy 
in an exterior algebra context. 

\vspace{1mm}
(2) Let $\M$ be a parameter space for objects in algebraic geometry 
such as the Hilbert scheme or a moduli space. 
Suppose that $\M$ is a subvariety of a rational variety $\GG$ of codimension $c$.   
Then the probability for a point $p$ in $\GG(\F_q)$ to lie in $\M(\F_q)$ is about $(1:q^c)$. 
This approach will be successful if the codimension $c$ is small and 
 the time required to check $ p \not\in \M(\F_q)$ 
is sufficiently small as compared with $q^c$.  
This technique was applied first by Schreyer~\cite{schreyer} 
to find  four different families of smooth surfaces in $\P^4$ 
with degree $11$ and sectional genus $11$ over $\F_3$ by a random search,   
and he provided a method to establish the existence of lifting these surfaces 
to characteristic zero. 
This technique has been successfully applied to solve various problems in constructive algebraic 
geometry (see \cite{st},  \cite{tonoli} and \cite{bel}). 

\vspace{1mm}
Here we illustrate the use of Techniques (1) and (2) in the first construction. 
For a fixed linear part of the second map $B$ 
in (\ref{eq:minimal-monad}), we find a map $B_2$ from $2\Omega^2(2)$ to $3\cO$ 
satisfying the conditions derived from the corresponding Tate resolution. 
These conditions gives an implicit algebraic description of the set $\M$  
of isomorphic classes of such $B_2$'s.  
In particular, we can show that $\M$ can be regarded 
as a four-codimensional subvariety $\M$ of  the grassmannian of $4$-quotient spaces 
$\GG$ of a ten-dimensional vector space.
This $B_2$ together with the linear part determines the Beilinson monad uniquely 
up to isomorphisms, and hence it suffices to find a point in $\M$ 
to construct a smooth surface.  
A point in $\M(\F_q)$ can be expected to be found 
from $\GG(\F_q)$ at a rate of $(1:q^4)$.
We use a random search over $\F_5$ to find a point in $\M(\F_5)$ 
and check smoothness of the corresponding surface $X/\F_5$. 
The existence of a lift to characteristic zero follows from the existence 
of a number field $\LL$ and a prime $\mathfrak{p}$ with 
$\cO_{\LL,\mathfrak{p}}/\mathfrak{p} \cO_{\LL,\mathfrak{p}}\simeq \F_5$. 
Then the generic point of $\Spec{\LL}$ corresponds to a smooth surface 
over $\LL$.  
 
The second construction is more deterministic.      
Let $A_1$ and $B_1$ be the linear parts of $A$ and $B$ respectively. 
In the example we found in the first construction, $A_1$ corresponds 
to a rational normal curve  in $\P(V)$; while $B_1$ corresponds 
to a rational cubic scroll in $\P(V)$, that does not intersect the rational 
normal curve associated to $A_1$. 
From these $A_1$ and $B_1$, one can reconstruct a smooth rational surface in $\P^4$
with the same invariants. Indeed,  the condition that  the composite of $B$ and $A$
is zero gives rise to the homogeneous system of $120$ linear equations 
with $140$ unknowns. In this case, the solution space has dimension 26. Thus,   
we can choose 26 variables freely to determine $A$ and $B$. 
A random choice of values for 26 parameters gives the Beilinson monad 
of type (\ref{eq:minimal-monad}), and the homology of this monad is 
the twisted ideal sheaf of the desired smooth surface.  
Let $\mathfrak{F}$ be the family of rational normal curves in $\P(V)$. 
For a fixed rational cubic scroll, and for each $ 20 \leq N \leq 26$, 
we can find the subfamily $\mathfrak{F}(N)$ of rational normal curves 
such that 
the associated homogeneous system of linear equations has the $N$-dimensional 
solution space.  The subfamily $\mathfrak{F}(N)$ is of codimension $N-20$ in $\mathfrak{F}$. 
So we can expect to find a point $p$ in $\mathfrak{F}(N)(\F_q)$ from $\mathfrak{F}(\F_q)$ 
at a rate of $\left(1:q^{N-20}\right)$.
The existence of a smooth rational surface in $\P^4$ with the desired invariants 
are established, when $N=22, 23,24,25$ and $26$.  

The calculations were done 
with the computer algebra system $\Mac$ 
developed by Grayson and Stillman~\cite{gs}. 
All the $\Mac$ scripts needed to construct surfaces are available at 
{\tt http://www.math.colostate.edu/$\sim$abo/programs.html}~\cite{abo}.    
  
\vspace{2mm}
\noindent 
{\it Acknowledgments.} 
I would like to thank Kristian Ranestad and Frank-Olaf Schreyer for 
very useful discussion about the subject of this article. 
The second construction was developed in collaboration with K. Ranestad;  
I received much help in working out this construction from him. 
The basic idea of the first construction was due to F.-O. Schreyer. 
He also  supported a great part of my computer skill development.   
Furthermore, I would like to thank them for permission to publish our results.

\section{Exterior algebra method}
Both Constructions I and II of rational surfaces use the technique of ``Beilinson monad". 
A Beilinson monad represents a given coherent sheaf 
in terms of direct sum of (suitably twisted) bundles of differentials 
and homomorphisms between these bundles, 
which are given by homogeneous matrices over an exterior algebra $E$. 
Recently, Eisenbud, Fl{\o}ystad and Schreyer~\cite{efs} showed 
that for a given sheaf, one can get the Beilinson monad from its ``Tate resolution",    
that is a free resolution over $E$, by a simple functor. 
This enables us to discuss the Beilinson monad in an exterior algebra  context. 
In this section, we take a quick look at the exterior algebra method  
developed by Eisenbud, Fl{\o}ystad and Schreyer. 
\subsection{Tate resolution of a sheaf}
Let $W$ be a $(n+1)$-dimensional vector space over a field $K$, 
let $V$ be its dual space, 
and let $\{x_i\}_{0 \leq i \leq n}$  and 
$\{e_i\}_{0 \leq i \leq n}$ be dual bases of $V$ and $W$ respectively. 
We denote by $S$ the symmetric algebra of $W$  and by $E$ 
the exterior algebra $\bigwedge V$ on $V$. 
Grading on $S$ and $E$ are introduced by $\deg(x)=1$ 
for $x \in W$ and $\deg(e)=-1$ for $e \in V$ respectively. 
The projective space of $1$-quotients of $W$ will be denoted by 
$\P^n=\P(W)$. 

Let $M=\bigoplus_{i\in \integer} M_i$ be a finite generated $S$-graded module. 
We set 
\[
 \omega_E:=\HHom_K(E,K)\simeq E \otimes_K \bigwedge^{n+1} W \simeq E(-n-1)
\]
and 
\[
 F^i:=\HHom_K(E,M_i)\simeq M_i \otimes_K \omega_E. 
\]
Let $\phi_i : F^i \rightarrow F^{i+1}$ be the map defined 
as a map taking $\alpha \in F^i$ to 
\[
\left[e \mapsto \sum_i x_i\alpha(e_i \wedge e) \right].  
\] 
Then the sequence 
\[
 \mathbf{R}(M) : \hspace{1mm} \cdots \rightarrow F^{i-1} \stackrel{\phi_{i-1}}{\lra} 
 F^i  \stackrel{\phi_i}{\longrightarrow} F^{i+1} \rightarrow \cdots 
\]
is a complex. This complex is eventually exact. Indeed,  
$\mathbf{R}(M)$ is exact at $\HHom_K(E,M_i)$ for all $i \geq s$ 
if and only if $s>r$,
where $r$ be the Castelnouvo-Mumford regularity of $M$ (see~\cite{efs} for a detailed proof).  
So starting from $\mathbf{T}(M)^{>r}:=\mathbf{T}(M_{> r})$, 
we can construct a doubly infinite exact $E$-free complex $\mathbf{T}(M)$ 
by adjoining a minimal free resolution of the kernel of $\phi_{r+1}$: 
\[
\mathbf{T}(M): \cdots 
\rightarrow T^r \rightarrow T^{r+1}:=\HHom_K(E,M_{r+1}) 
\stackrel{\phi_{r+1}}{\longrightarrow} \HHom_K(E,M_{r+2}) \rightarrow \cdots.  
\]
This $E$-free complex is called the {\it Tate resolution} of $M$. 
Since $\mathbf{T}(M)$ can be constructed by starting from $\mathbf{R}(M_{>s})$, 
$s\geq r$, the Tate resolution depends only on the sheaf 
$\sF = \widetilde{M}$ on $\P(W)$ associated to $M$.
 We call $\mathbf{T}(\sF):=\mathbf{T}(M)$ the {\it Tate resolution} of $\sF$. 
The following theorem gives a description of all the term of a Tate resolution: 
\begin{theorem}[\cite{efs}]
Let $M$ be a finitely generated graded $S$-module and 
let $\sF := \widetilde{M}$ be the associated sheaf on $\P(W)$.  
Then the term of the complex $\mathbf{T}(\sF)$ with cohomological degree $i$ 
is $\bigoplus_j \HH^j \sF(i-j) \otimes \omega_E$.
\label{th:tate-resolution}
\end{theorem}
\subsection{Beilinson monad}
Eisenbud, Fl{\o}ystad and Schreyer~\cite{efs} showed that 
applying a simple functor to the Tate resolution $\mathbf{T}(\sF)$ 
gives  a finite complex of sheaves whose homology is 
the sheaf $\sF$ itself:  
Given $\mathbf{T}(\sF)$, we define $\Omega(\sF)$ to be the complex 
of vector bundles on $\P(W)$ obtained by replacing each 
summand $\omega_E(i)$ by the $i$th twisted bundle $\Omega^i(i)$ of differentials. 
The differentials of the complex are given by using isomorphisms 
\[
 \HHom_E(\omega_E(i),\omega_E(j)) \simeq \bigwedge^{i-j} V 
 \simeq \HHom(\Omega^{i}(i),\Omega^j(j)).
\]
\begin{theorem}[\cite{efs}]
 Let $\sF$ be a coherent sheaf on $\P(W)$. 
 Then $\sF$ is the homology of $\Omega(\sF)$ in cohomological 
 degree $0$, and $\Omega(\sF)$ has no homology otherwise. 
\label{th:beilinson-monad}
\end{theorem}

\noindent 
We call $\Omega(\sF)$ the {\it Beilinson monad} for $\sF$. 
\section{Construction I }
This section will be devoted to constructing 
a family of  rational surfaces $X$ in $\P^4$ with degree $d=12$,  
sectional genus $\pi=13$ and $p_g=q=0$.   
The construction discussed in this section takes the following  
three steps:
\begin{itemize}
\item[(1)] 
Find a smooth surface $X$ with the prescribed invariants 
over a finite field of a small characteristic. 
\item[(2)] 
Determine the type of the linear system, 
which embeds $X$ into $\P^4$  to justify that  
the surface $X$ found in the previous step is rational.  
\item[(3)] 
Establish the existence of a lift to characteristic zero. 
\end{itemize}

\subsection{Construction over a small field}
Let $K$ be a field, let $W$ be a five-dimensional vector space over $K$ 
with basis $\{x_i\}_{0 \leq i \leq 4}$, 
and let $V$ be its dual space with dual basis $\{e_i\}_{0 \leq i \leq 4}$. 
Let $X$ be a smooth surface in $\P^4=\P(W)$ with the invariants given above. 
The first step is to determine the type of the Beilinson monad for the twisted 
ideal sheaf of $X$, 
which is derived from the partial knowledge of its cohomology groups. 
Such information can be determined from 
general results such as the Riemann-Roch formula and Kodaira vanishing theorem 
(see~\cite{des} for more detail).  
Here we assume that $X$ has the natural cohomology in the range $-1 \leq j \leq 4$   
of twists:
{
$$
{
\setlength{\br}{10mm} 
\setlength{\ho}{6mm} 
\fontsize{10pt}{11pt} 
\selectfont 
\begin{xy} 
%
%
<-3.5\br,0\ho>;<5.5\br,0\ho>**@{-}?>*@{>}?(0.9)*!/^3mm/{j} 
%
%
,<-2\br,0\ho>;<-2\br,6\ho>**@{-}?>*@{>}?(0.9)*!/^3mm/{i} 
%
%
,0+<-3.5\br,5\ho>;<4.5\br,5\ho>**@{-}
,0+<-3\br,4\ho>;<4\br,4\ho>**@{-}
,0+<-3\br,3\ho>;<4\br,3\ho>**@{-}
,0+<-3\br,2\ho>;<4\br,2\ho>**@{-}
,0+<-3\br,\ho>;<4\br,\ho>**@{-}
%
%
,0+<-3\br,0\ho>;<-3\br,5\ho>**@{-} 
,0+<-2\br,0\ho>;<-2\br,5\ho>**@{-} 
,0+<-1\br,0\ho>;<-1\br,5\ho>**@{-} 
,0+<3\br,0\ho>;<3\br,5\ho>**@{-} 
,0+<2\br,0\ho>;<2\br,5\ho>**@{-} 
,0+<1\br,0\ho>;<1\br,5\ho>**@{-}   
,0+<0\br,0\ho>;<0\br,5\ho>**@{-} 
,0+<4\br,0\ho>;<4\br,5\ho>**@{-} 
%
%
,0+<-2.5\br,3.5\ho>*{13}
,0+<-.5\br,2.5\ho>*{4}
,0+<.5\br,2.5\ho>*{2}
,0+<1.5\br,1.5\ho>*{2}
,0+<2.5\br,1.5\ho>*{3}
,0+<3.5\br,0.5\ho>*{5}
,0+<5.5\br,2.5\ho>*{\h^i \sI_X(j)}
%
%
,<-2.5\br,-.5\ho>*{-1}
,<-1.5\br,-.5\ho>*{0}
,<-.5\br,-.5\ho>*{1}
,<.5\br,-.5\ho>*{2}
,<1.5\br,-.5\ho>*{3}
,<2.5\br,-.5\ho>*{4}
,<3.5\br,-.5\ho>*{5}
%
%
,<4.2\br,.5\ho>*{0}
,<4.2\br,1.5\ho>*{1}
,<4.2\br,2.5\ho>*{2}
,<4.2\br,3.5\ho>*{3}
,<4.2\br,4.5\ho>*{4}
\end{xy}
}
$$
}

\noindent
Here a zero is represented by the empty box. 
By Theorem~\ref{th:tate-resolution},
the Tate resolution $\mathbf{T}(\sI_X)[4]=\mathbf{T}(\sI_X(4))$ 
includes an exact $E$-free complex of the following type:
\begin{eqnarray}
\cdots \rightarrow 13\omega_E(5) 
\rightarrow 4\omega_E(3) \rightarrow 2\omega_E(2)\oplus 2\omega_E(1) 
\rightarrow 3\omega_E \oplus 5 \omega_E(-1) \rightarrow \cdots. 
\label{eq:betti-diagram}
\end{eqnarray}
\noindent
From Theorem~\ref{th:beilinson-monad}, 
it follows, therefore, that   the corresponding Beilinson monad for 
$\sI_X(4)$ is of the following type: 
\begin{eqnarray}
0 \rightarrow 4\Omega^3(3) \stackrel{A}{\rightarrow}
 2\Omega^2(2)\oplus2\Omega^1(1) 
\stackrel{B}{\rightarrow} 3\cO \rightarrow 0. 
\label{eq:monad}
\end{eqnarray}

\vspace{2mm}
The next step is to describe what maps $A$ and $B$ 
could be the differentials of  Monad~(\ref{eq:monad}). 
The identifications  
\[
\HHom(\Omega^i(i),\Omega^j(j)) \simeq \HHom_E(\omega_E(i), \omega_E(j)) \simeq \HHom_E(E(i),E(j)), 
\]
allow us to think of the maps $A$ and $B$ as homomorphisms 
between $E$-free modules.  
Since the composite of $B$ and $A$ is zero, 
each column of $A$ can be written as an $E$-linear combination 
of columns of $\Syz(B)$. 
A reasonable case is when the minimal free resolution $\mathbf{B}$ of $\Coker(B)$ 
has the following type: 
\begin{eqnarray}
 \begin{tabular}{c|ccccc}
 \hline
 0 & 3 & 2 & . & .\\
 -1 & . & 2 & 4  & $a_2$  \\
 -2 & . & . & $a_1$ & $\ast$
 \end{tabular} 
 \label{eq:betti-table}
\end{eqnarray}
Assume that there exists such a map $B$. 
Then this $B$ uniquely determines $A_B : 4E(3) \rightarrow 2E(2)\oplus 2E(1)$, 
which could be the first map in (\ref{eq:monad}), 
up to automorphism of  $4E(3)$. 
In order to exclude superfluous candidates, 
we check conditions the map $B$ and the associated map $A_B$ should satisfy. 
\begin{proposition}
Suppose that $\mathbf{B}$ is minimal. Then 
 $a_1=5$ and $a_2=0$. 
\label{th:type-of-syzygies}
\end{proposition} 
\begin{proof}
The exact sequence $\mathbf{B}[1]$ includes 
the linear transformation 
\[
B(1):  \HH^0 2\Omega^2(3) \oplus \HH^0 2\Omega^1(2) \rightarrow 
 \HH^0 3\cO(1). 
\]
Recall that $\HH^0 \Omega^i(i+1)$ is isomorphic to $\bigwedge^{i+1}W$ 
and there are canonical identifications 
$\bigwedge^{i+1} W  \simeq \bigwedge^{4-i} V$. 
From (\ref{eq:betti-table}), we can deduce that $\Ker(B(1))$ is isomorphic to 
the quotient of $4V \oplus  a_1 K$ by $a_2 K$. 
Since $\mathbf{B}$ is minimal, 
the map $a_2 E(4) \rightarrow a_1 E(3)$ in (\ref{eq:betti-table}), 
and hence the corresponding linear transformation 
$a_2 K \rightarrow a_1 K$,  is zero. 
Thus the quotient is $4V/a_2K \oplus  a_1K$. 
However this is possible only if $a_2=0$, 
because $A_B(1)$ induces an injective linear transformation 
from $\HH^0 4\Omega^3(4)$ to $\Ker(B(1))$. 
By assumption,  we obtain
\begin{eqnarray}
\h^0 \sI_X(5)=\dim(\Ker(B(1)))-4\cdot \h^0 \Omega^{3}(3)=(20+a_1)-20=a_1,    
\label{eq:ker}
\end{eqnarray}
and thus $a_1=5$. 
\end{proof}
Consider the linear transformation 
\[
A_B(-5) : \HH^4 4\Omega^3(-2) \rightarrow \HH^4 2\Omega^2(-3) \oplus \HH^4 2\Omega^1(-4).  
\]
The kernel of this linear transformation is isomorphic to $\HH^3 \sI_X(-1)$. 
From the fact that $(\Omega^i)^{\vee} \simeq \Omega^{4-i}(5)$ 
and Serre duality it follows that for all $i=1,2,3$, there are isomorphisms 
\[
 (\HH^4 \Omega^i(i-5))^* \simeq \HH^0 \Omega^{4-i}(5-i) \simeq \bigwedge^i V. 
\]
Hence $\HH^3 \sI_X(-1)$ has the form $b_1 V \oplus b_2  K$. 
By assumption, the syzygy matrix of $A_B$ does not contain 
any linear entries. 
This implies that the cokernel of $A_B$ 
has the minimal free presentation of the following shape: 
\begin{eqnarray}
 \begin{tabular}{c|cccc}
 \hline
 -1 & 2 & . & . \\
 -2 & 2& 4 & .    \\
 -3 & . & . & 13 
  \end{tabular} 
 \label{eq:betti-table-2}
\end{eqnarray}
After all, our task is to find maps $B$ satisfying Conditions (\ref{eq:betti-table}), 
(\ref{eq:ker}) such that $A_B$ fulfills Condition (\ref{eq:betti-table-2}). 

\vspace{2mm}
Let $B_1$ be the linear part of $B$. 
Assume that $B_1$ is the general member of $\HHom(2E(1),3E)$. 
Consider the following problem:

\vspace{2mm}
\noindent 
{\bf Problem: } 
Find a $B_2 \in \HHom(2E(2),3E)$ such that 
$B=(B_2,B_1)$ satisfies Condition (\ref{eq:betti-table}). 

\vspace{2mm}
\noindent  
To ease our calculations, 
we define $B_1$ by the matrix 
\[
B_1=
\left(
\begin{array}{cc}
e_0 & e_1 \\
e_1 & e_2 \\
e_3 & e_4 
\end{array}
\right), 
\] 
and for this $B_1$, we detect a $B_2\in \HHom_E(2E(2),3E)$ 
satisfying (\ref{eq:betti-table}).   
We then check that such a $B$ actually gives rise to a smooth surface $X$  
and the surface is rational. 
The type of a linear system embedding $X$ into $\P^4$ 
will be also determined.

\begin{remark}
(i)  Note that $B_1$ corresponds to the multiplication map 
\[
 W \tensor \HH^1 \sI_X(3) \rightarrow \HH^1\sI_X(4)
\]
(see \cite{efs}). 
So $B_1$ induces a map from the set of hyperplanes 
to the set of linear transformations from $\HH^1 \sI_X(3)$ to $\HH^1 \sI_X(4)$. 
By the definition of $B_1$, the hyperplane classes,  
for which the corresponding linear transformation is not injective, 
form a rational cubic scroll in $\P(V)$.

\vspace{1mm} 
\noindent
(ii)  
Let   
\[
 B_1=
\left(
\begin{array}{cc}
e_0 & e_1 \\
e_1 & e_2 \\
e_2 & e_3 
\end{array}
\right).  
\]
The associate variety in $\P(V)$ is a cone over 
a twisted cubic, so this is singular. 
For this $B_1$, 
there are actually $B_2$'s such that $(B_2,B_1)$ give rise to 
monads of the desired type by using Construction II we will discuss later.    
In these examples, the corresponding surfaces have, however, singularities. 
So our choice of the matrix $B_2$ seems to be reasonable. 
\end{remark}
Both the columns of such a $B_2$ are not generated by those of $B_1$.  
They are furthermore linearly independent.   
So $B_2$ corresponds to a plane in the quotient of 
$3 \bigwedge^2 V$ by $V \wedge <B_1> $.  
Obviously, 
$B_2$ and $B_2'$ in $\HHom(2E(2),3E)$ corresponds to the same plane 
in the quotient if and only if 
$B=(B_2,B_1)$ and $B'=(B_2',B_1)$ define the same $E$-module. 
In this case, $(A_B,B)$ and $(A_{B'},B')$ give rise to the same monad 
up to  an isomorphism. 
Since the quotient space has dimension $20$, 
elements in $\HHom(2E(2),3E)$ satisfying the desired conditions 
form a subvariety of $G(2,20)$.  

Let us try to describe the variety $\M$ more precisely.  
From (\ref{eq:betti-table}), the following sequence of vector spaces 
is exact:
\begin{eqnarray}
\begin{CD}
 0 \rightarrow 4K \rightarrow 2  V \oplus 
 2\bigwedge^2 V \stackrel{B}{\rightarrow} 3 \bigwedge^3 V \stackrel{f}{\rightarrow} U \rightarrow 0, 
\end{CD}
\label{eq:exact-squence-vector}
\end{eqnarray}
where $U$ is the cokernel $\Coker(B_2,B_1)$ that is four-dimensional. 
Note that the map $B_2$ is factored as in the following figure: 
$$
\xymatrix{
& 2K \tensor V  
\ar[dr]_{B_2 \tensor 1} 
\ar[rr]^{B_2}    
&&  3 \bigwedge^3 V. \\
&& 3 \bigwedge^2 V  \tensor V \ar[ur]_{\wedge} \\
}
$$
Define a map $g$ from $3 \bigwedge^2 V \tensor  V$ to $U$ as 
the composite of  the map $f$ from $3\bigwedge^3 V$  to $U$ in 
(\ref{eq:exact-squence-vector}) and the map 
$\wedge : 3 \bigwedge^2 V \tensor V \rightarrow 3\bigwedge^3 V$ 
in the above figure.  
This map corresponds to a map $g^* \in \HHom(3 \bigwedge^2 V, U \otimes W)$, and   
each column of $B_2$ is obtained as 
an element of $\Ker(g)$. 
By construction, $\dim(\Ker(g^*))=12$, and hence $\dim(\M)=\dim G(2,12)=20$. 
So $\mathrm{codim}(\M,G(2,20))=36-20=16$. 
Fortunately, we can improve our description 
of $\M$ to obtain smaller codimension as follows: 
Let $T$ be the cokernel $\Coker(A_1)$. 
Then one sees that $T$ is isomorphic to $10 \bigwedge^5 V$. 
We denote the grassmaniann of $4$-quotient spaces of $T$ by $\GG$.  
\begin{proposition}
$\M$ is a four-codimensional subvariety of $\GG$.
\end{proposition}
\begin{proof}
Consider the following commutative diagram:
\[
\begin{CD}
& & 0 @>>> 2 \bigwedge^2 V @>B_1>>  3\bigwedge^3 V @>>> T @>>> 0 \\
& & & & @VVV  \| &  & @VVV \\
0 @>>> 4 K @>>> 2V \oplus 2\bigwedge^2V @>(B_2,B_1)>> 
3 \bigwedge^2 V @>>> U @>>> 0
\end{CD}
\]
So  there is a map from $\HHom_K(2V, 3\bigwedge^3V)$ to the set 
of $4$-quotient spaces $\GG$ of $T$ defined by $B_2 \mapsto U$. 
This map induces an injective map from $\M$ into $\GG$. 
Indeed, $B_2$ and $B_2'$ are mapped to the same element of $\GG$ 
if and only if  the cokernels of $(B_2,B_1)$ and $(B_2',B_1)$ are isomorphic. 
This happens precisely if 
each column of $B_2'$ is an $E$-linear combination of 
the columns of $B_1$ and $B_2$. 
Therefore $\M$ can be regarded as a four-codimensional subvariety of $\GG$.  
\end{proof}

For a field $\F_q$ with $q$ elements, we can therefore expect to find a point 
$\M(\F_q)$ from $\GG(\F_q)$ at a rate of $(1:q^4)$. 
The statistics suggests that there two different family in 
$\GG(\F_5)$ whose elements satisfy Condition (\ref{eq:betti-table}):  
\[
 \begin{tabular}{c|ccccc}
 \hline
 2& 4 &. &. &. \\
 1 &. & 3 & 2 & .\\
 0 &. & . & 2 & 4   \\
 -1 &. & . & . & 5 
 \end{tabular} 
 \quad \mbox{and} \quad 
  \begin{tabular}{c|ccccc}
 \hline
 2& 4 &. &. &. \\
 1 &. & 3 & 2 & .\\
 0 &. & . & 2 & 4   \\
 -1 &. & . & . & 10 
 \end{tabular}  
\]
However, examples with the second Betti table do not give rise to 
any Beilinson monads of the desired type by Proposition~\ref{th:type-of-syzygies}.  
Examples with the first Betti table appeared $18$ times in a test of $5^4\cdot 10$ examples. 
So it turns out that the codimension of the first family in $\GG(\F_5)$ is $4$.

\begin{proposition}
There is a smooth surface in $\P^4$ over $\overline{\F}_5$ with
$d=12$ and $\pi=13$. 
\label{th:small-field}
\end{proposition}
\begin{proof}
By random search, 
we can find a $B_2$ satisfying the desired conditions. 
Determine the corresponding maps 
$A_B: 4\Omega^3(3) \rightarrow 2\Omega^2(2) \oplus 2\Omega^1(1)$ 
and $B=(B_2,B_1) : 2\Omega^2(2) \oplus 2\Omega^1(1) \rightarrow 3\cO$. 
Then compute the homology $\ker(B)/\im(A_B)$. 
If the homology is isomorphic to the ideal sheaf of a surface 
with the desired invariants, then check smoothness 
of the surface with the Jacobian criterion. 
Perform this operation repeatedly.  
Then one can find the matrix 
\[
 B_2=
 \left(
 \begin{array}{cc}
 e_{23}-e_{34} & 2e_{23}+e_{24}-2e_{34} \\
 e_{03} + e_{13} -2e_{14} + e_{24} - 2e_{34} & 
 2e_{13} -2 e_{23} +2e_{04} + e_{24} + e_{34}   \\
 b_{31} & b_{32} 
 \end{array}
 \right), 
\]
where 
\begin{eqnarray*} 
  b_{31} &=& e_{01} -2e_{02} -2e_{12} + 2e_{03} - 
                      2e_{13}+e_{04} -2e_{14} -e_{24} 
                      -  2e_{34}; \\              
 b_{32} &=& e_{02} -e_{12} -e_{03} + e_{23} - 
                      e_{04}+2e_{14}  - 2e_{34}.                                          
\end{eqnarray*}
Here we denote $e_i \wedge e_j$ by $e_{ij}$ for each $0\leq i<j \leq 4$. 
A script for finding $B_2$ can be obtained from our webpage~\cite{abo}. 
For this $B_2$, there is a smooth surface 
with $d=12$ and $\pi=13$. 
\label{th:first-construction}
\end{proof}
\subsection{Adjunction process}
\label{sec:adjunction_process}
In this subsection, we spot the surface found in the previous step   
within the Enriques-Kodaira classification 
and determine the type of the linear system that embeds $X$ into $\P^4$. 
First of all, 
we recall a result of Sommese and Van de Ven for a surface over $\C$: 
\begin{theorem}[\cite{sv}]
Let $X$ be a smooth surface in $\P^n$ over $\C$ 
with degree $d$, sectional genus  $\pi$, geometric genus $p_g$ and irregularity $q$,  
let $H$ be its hyperplane class, 
let $K$ be its canonical divisor and let $N=\pi - 1 + p_g -1$. 
Then the adjoint linear system $| H+K |$ defines a birational morphism 
\[
 \Phi=\Phi_{|H+K|} : X \rightarrow \P^{N-1}
\]
onto a smooth surface $X_1$, which blows down precisely all $(-1)$-curves on $X$, 
unless 
\begin{itemize}
 \item[$(i)$] $X$ is a plane, or Veronese surface of degree $4$, or $X$ is ruled by lines; 
 \item[$(ii)$] $X$ is a Del Pezzo surface or a conic bundle;
 \item[$(iii)$] $X$ belongs to one of the following four families: 
 \begin{itemize} 
  \item[$(a)$]  $X=\P^2(p_1, \dots , p_7)$ embedded by $H \equiv 6L-\sum_{i=0}^7 2E_i$;
  \item[$(b)$]  $X=\P^2(p_1, \dots , p_8)$ embedded by $H \equiv 6L-\sum_{i=0}^7 2E_i-E_8$;
  \item[$(c)$]  $X=\P^2(p_1, \dots , p_8)$ embedded by $H \equiv 9L-\sum_{i=0}^8 3E_i$;
  \item[$(d)$] $X=\P(\sE)$, where $\sE$ is an indecomposable rank $2$ bundle over an
  elliptic curve and $H\equiv B$, where $B$ is a section $B^2=1$ on $X$. 
 \end{itemize}
\end{itemize}
\label{th:adjoint-linear-system}
\end{theorem}
\begin{proof}
 See~\cite{sv} for the proof. 
\end{proof}

\noindent 
Setting $X=X_1$ and performing the same operation repeatedly, 
we obtain a sequence 
\[
 X \rightarrow X_1 \rightarrow X_2 \rightarrow \cdots \rightarrow X_k. 
\] 
This process will be terminated if $N-1 \leq 0$. 
For a surface with nonnegative Kodaira dimension, 
one obtains the minimal model at the end of the adjunction process. 
If the Kodaira dimension equals $-\infty$, 
we end up with a ruled surface, 
a conic bundle, a Del Pezzo surface, 
$\P^2$, or one of the few exceptions of Sommese and Van de Ven. 

It is not known whether the adjunction theory holds over a finite field. 
However, we have the following proposition:
\begin{proposition}[\cite{ds}]
Let $X$ be a smooth surface over a field of arbitrary  characteristic. 
Suppose that the adjoint linear system $|H+K|$ is base point free. 
If the image $X_1$ in $\P^N$ under the adjunction map $\Phi_{|H+K|}$ 
is a surface of the expected degree $(H+K)^2$, the expected sectional genus 
$\frac{1}{2}(H+K)(H+2K)+1$ and  with $\chi(\cO_X)=\chi(\cO_{X_1})$, 
then $X_1$ is smooth and $\Phi: X\rightarrow X_1$ is a simultaneous 
blow down of the $K_1^2-K^2$ many exceptional lines on $X$. 
\label{th:adjunction-p}
\end{proposition}
\begin{proof}
See Proposition 8.3 in~\cite{ds} for a proof. 
\end{proof}
\begin{remark}
The exceptional divisors contracted in each step is defined over the base field.
\end{remark}
In~\cite{des} and~\cite{ds}, it is described 
how to compute the adjunction process for a smooth surface given by explicit equations 
(see~\cite{ds} for the computational details).   
Let $X$ be the smooth surface found in the previous step. 
The computation for the adjunction process in characteristic $5$ gives 
\begin{eqnarray}
H\equiv 12L-\sum_{i_1=1}^2 4E_{i_1}-\sum_{i_2=3}^{11} 3E_{i_2} 
-\sum_{i_3=12}^{14} 2E_{i_3}-\sum_{i_4=15}^{21} E_{i_4}, 
\label{eq:ample-divisor}
\end{eqnarray}
where $L$ is the class of a line in $\P^2$. 
This process ends up with the Del Pezzo surface of degree $7$, 
which is a two points blowing up of $\P^2$. 
Therefore we can conclude that $X$ is rational. 

\begin{remark}
Another way to prove rationality of $X$ is to count the number of 
$6$-secant lines to $X$.
First we prove the following claim: 
Let $Y$ be a smooth surface in $\P^4$ with $d=12$, $\pi=13$ and 
$p_g=q=0$. Suppose that there are at most a finite number of $6$-secant 
lines to $Y$. Then $Y$ is rational. 

To prove this claim, we use 
the Le Barz's formula $N_6(d,\pi, \chi)$. This formula gives us $N_6(12,13,1)=8$,  
that equals the number 
of $6$-secant lines to $Y$ plus the number of exceptional lines 
on $Y$, if there are at most a finite number of $6$-secant lines to $Y$ and 
if there are no lines on $Y$ with nonnegative self-intersection number~\cite{leBarz}. 

Let $H_2$ be the hyperplane class of the second adjoint surface 
$Y_2$ and let $K_2$ be the canonical divisor of $Y_2$. 
Then we have $H_2 \cdot K_2= -12+a$, 
where $a$ is the number of exceptional lines on $Y$. 
The Le Barz's formula tells us that $Y$ can have at most eight exceptional lines 
if there are a finite number of $6$-secant lines to $Y$. 
So $H_2 \cdot K_2 < 0$, and hence $Y$ is rational.

Next we show that $X$ has only one $6$-secant line. 
The union of $6$-secant lines to $X$ is contained in 
all the quintics that contain $X$. 
With $\Mac$, we can check that 
\[
 V(\HH^0 \sI_X(5))=X \cup L_0, 
\]
where $L_0$ is a line. 
So $L_0$ is the only $6$-secant line to $X$, 
and hence $X$ is rational. 
\end{remark}

\subsection{Lift to characteristic zero} 
\label{lift}
In the previous step, we constructed a smooth surface in $\P^4$ over $\bar{\F}_5$. 
However, our main interest is the field of complex numbers $\C$.  
In this section,  
we show the existence of a lift to characteristic zero as follows: 
Let $\M$ and $\GG$ be given as in the previous subsection.  
\begin{proposition}[\cite{schreyer}]
Let $y \in \M(\F_q)$ be a point, 
where $\M(\F_q)$ has codimension $4$. 
Then there exists  a number field $\LL$ and  
a prime $\mathfrak{p}$ in $\LL$  such that the residue field 
$\cO_{\LL,\mathfrak{p}} /\mathfrak{p}\cO_{\LL,\mathfrak{p}}$ is $\F_p$.
Furthermore, if the surface $X/\F_p$ corresponding to $y$ is smooth, 
then the surface $X/\LL$ corresponding to the generic point 
of $\Spec \LL \subset \Spec \cO_{\LL,\mathfrak{p}}$ is also smooth. 
\label{th:lift}
\end{proposition}
\begin{proof}
Let $q$ be a prime number. 
If this is not the case, $\Z$ has to be replaced by the ring of integers 
in a number field which has $\F_q$ as the residue field. 

Since $\M(\F_q)$ has pure codimension $4$ in $y$, 
there are four hyperplanes $H_1,\dots,H_4$ in $\GG$ 
such that $y$ is an isolated point of $\M(\F_q) \cap H_1 \cap \cdots \cap H_4$. 
We may assume that $H_1, \dots, H_4$ are defined over $\Spec\Z$ 
and that they meet transversally in $y$. 
This allows us to think that $\M(\F_p) \cap H_1 \cap \cdots \cap H_4$ 
is defined over $\Z$. Let $Z$ be an irreducible component of $\M_{\Z}$.
Then $\dim Z=1$. 

The residue class field of generic point of $Z$ is a number field $\LL$ that is 
finitely generated over $\Q$, 
because $\M_{\Z}$ is projective. 
Let $\cO_{\LL}$ be the ring of integers of $\LL$ 
and let $\mathfrak{p}$ be a prime ideal which lies over $y \in Z$. 
Then $\Spec\cO_{\LL, \mathfrak{p}} \rightarrow Z \subset \M$ 
is an $\cO_{\LL,\mathfrak{p}}$-valued point which lifts $y$. 
The residue class field $\cO_{\LL,\mathfrak{p}}/\mathfrak{p}\cO_{\LL,\mathfrak{p}}$ 
is a finite extension of $\F_p$. 

Performing the construction of the surface over $\cO_{\LL,\mathfrak{p}}$ 
gives a flat family $\sX$ of surfaces over $\cO_{\LL,\mathfrak{p}}$. 
Since smoothness is an open property, and since the special fiber 
$\sX_{\mathfrak{p}}$ is smooth, 
the general fiber $\sX_{\LL}$ is also smooth. 
\end{proof}
Next we argue that the adjunction process of the surface over the number field 
$\LL$ has the same numerical behavior: 
\begin{proposition}[\cite{ds}]
Let $\sX \rightarrow \Spec \cO_{\LL,\mathfrak{p}}$ be a family as in 
Proposition~\ref{th:existence-of-lift}. 
If the Hilbert polynomial of the first adjoint surface of $X=\sX \tensor \F_q$ 
is as expected, and if 
$\HH^1(X,\cO_X(-1))=0$, then 
the adjunction map of the general fiber $\sX_{\LL}$ blows down 
the same number of exceptional lines as the adjunction map of the spacial fiber $X$. 
\label{th:adjunction-0}
\end{proposition}
\begin{proof}
See Corollary 8.4 in~\cite{ds} for a proof. 
\end{proof}
\begin{theorem}
There exists a family of  smooth rational surfaces in $\P^4$ over $\C$ 
with $d=12$, $\pi=13$, $p_g=q=0$ and 
\[
H \equiv 12L-\sum_{i_1=1}^2 4E_i-\sum_{i_2=3}^{11} 3E_i 
-\sum_{i_3=12}^{14} 2E_i-\sum_{i_4=15}^{21} E_i. 
\]
\label{th:existence-of-lift}
\end{theorem}
\begin{proof}
Let $\M_{\F_5}=\M(\F_5)$ and let $y$ be the element of $\M_{\F_5}$ 
corresponding to the surface obtained in (1).   
We check with $\Mac$ that $y$ satisfies the condition in the previous proposition. 
Let 
\begin{eqnarray*}
V_1&=&\HHom(4\Omega^3(3),2\Omega^2(2) \oplus 2\Omega^1(1))
\simeq 8V \oplus 8\bigwedge^2 V, \\
V_2&=&\HHom(2\Omega^2(2) \oplus 2\Omega^1(1),3\cO) 
\simeq 6\bigwedge^2 V \oplus 6 V,   \\ 
V_3&=&\HHom(4\Omega^3(3),3\cO) \simeq 12 \bigwedge^3 V.
\end{eqnarray*}
Consider the map  
\[
\phi:  V_1 \oplus V_2 \rightarrow V_3 
\]
defined by $(A',B') \mapsto B' \circ A'$.  Let $A_B\in V_1$ and $B  \in V_2$ 
be the differentials of the monad for $X$ given in (1). 
The differential 
$d\phi: V_1 \oplus V_2 \rightarrow V_3$  of the map $\phi$  
at the point $(A_B,B)$ is given by $(A',B') \mapsto B \circ A'+B' \circ A_B$. 
Consider the subset  $\tilde{P}$ of $\phi(0)$ whose elements give 
monads of type (\ref{eq:monad}). This forms an open subset of $\phi(0)$.  
Let $H$ denote the group 
\[
 \left\{\left.
 \left( 
 \begin{array}{cc}
 C & 0  \\
 v &  D 
 \end{array}
 \right) 
 \right|
 C,D \in GL(2,\F_5),  v \in GL(2,V) \right\}
\]
and let $G=GL(4,\F_5) \times H \times GL(3,\F_5)$. 
Then $G'=G/\F_5^{\times}$ acts on $\tilde{P}$ by 
\[
(A',B')(f,g,h)=\left(g \circ A' \circ f^{-1},h \circ B' \circ g^{-1}\right). 
\]  
Let  $P$ be the set of isomorphic classes of monads of type (\ref{eq:monad}) 
and let $T_{P,(A_B,B)}$ be the Zariski tangent space 
of $P$ at the point corresponding to $(A_B,B)$.  
Then $P \simeq \tilde{P}/G'$, 
and hence $\dim \left(T_{P,(A_B,B)}\right)= \dim\left((d\phi)^{-1}(0)/G' \right)$. 
For fixed bases of $V_1, V_2$ and $V_3$, 
we can give the matrix that represents the differential $d\phi$ explicitly. 
This matrix enables us to compute the kernel of $d\phi$. 
This computation can be done with $\Mac$:  
\[
 \dim\left(T_{P,(A_B,B)}\right)=\dim\left((d\phi)^{-1}(0)\right)-\dim(G')=90-(53-1)=38. 
\]
A $\Mac$script for this computation can be found in~\cite{abo}.  
By construction, Monad (\ref{eq:monad}) is uniquely determined by $B$ 
up to isomorphisms. 
As $B_1$ corresponds to the general member of $\HHom(2E(-1),3E)$, 
the parameter space for $B_1$ has dimension $18$, 
which implies that $\M_{\F_5}$ has the desired dimension $38-18=20$. 

The type of the very ample divisor that embeds the surface into $\P^4$ 
follows from Proposition~\ref{th:adjunction-0}. 
\end{proof}
\section{Construction II} 
First of all, we motivate the second construction. 
Let $X$ be the smooth rational surface in $\P^4$ constructed in the previous section. 
Recall that the twisted ideal sheaf $\sI_X(4)$ of $X$ is obtained from the monad 
\[
 4\Omega^3(3) \stackrel{A}{\rightarrow} 2\Omega^2(2) \oplus 2\Omega^1(1) 
 \stackrel{B}{\rightarrow} 3\cO.  
\]
Let $A_1$ and $B_1$ be the linear parts of $A$ and $B$ respectively. 
The variety in $\P(V)$ associated to $A_1$ is a rational normal curve $C_{A_1}$;   
while the variety in $\P(V)$ associated to $B_1$ is a rational cubic scroll $X_{B_1}$. 
With $\Mac$, we can check that $X_{B_1}$ and $C_{A_1}$ 
do not intersect. 
Starting with these $A_1$ and $B_1$, we reobtain the surface $X$.  
The construction takes the following five steps $(1), (2), (3),(4)$ and $(5)$:

\vspace{2mm}
\noindent 
(1)  Detect an $A_2' \in \HHom(4\Omega^3(3),2\Omega^2(2))$
and a $B_2' \in \HHom(2\Omega^2(2),3\cO)$ 
such that the corresponding sequence 
\begin{eqnarray}
 4\Omega^3(3) \stackrel{A'}{\rightarrow} 2\Omega^2(2) \oplus 2\Omega^1(1) 
 \stackrel{B'}{\rightarrow} 3\cO  
\label{eq:sequence}
\end{eqnarray}
is a monad, 
where $A'={}^t(A_1,A_2')$ and $B'=(B_2',B_1)$. 
Since Sequence~(\ref{eq:sequence}) should be a complex,   
$A_2'$ and $B_2'$ must satisfy the following condition:  
\begin{eqnarray}
 B_2' \circ A_1 + B_1 \circ A_2'=0. 
\label{eq:relations}
\end{eqnarray}
Let $A_2'$ be the $2 \times 4$ matrix whose $(k,l)$ entry is 
$
 \sum_{i<j} a_{ij}^{kl}e_i\wedge e_j
$
and let $B_2'$ be the $3 \times 2$ matrix whose $(k,l)$ entry is 
$
 \sum_{i<j} b_{ij}^{kl}e_i\wedge e_j. 
$
Condition (\ref{eq:relations}) gives rise to the homogeneous 
system of $120$ linear equations with $140$ unknowns. 
The minimal number of equations in the system, 
denoted by $N_{A_1}$,  is $114$.  
Thus the solution space to the system has dimension $26$. 

\vspace{2mm}
\noindent
(2) Solving those equations for $b_{ij}^{kl}$'s, 
we obtain $34$ relations among $b_{ij}^{kl}$'s, and hence  
$26$ variables in $b_{ij}^{kl}$'s can be chosen freely to determine $B_2'$. 
With $\Mac$, 
it can be checked that $B_1$ and the $B_2'$ given by the random choices of 
values for $26$ parameters 
define a homomorphism from $2E(2)\oplus 2E(1)$ to $3E$ satisfying 
Conditions  (\ref{eq:betti-table}), (\ref{eq:betti-table-2})  
and the conditions in Proposition~\ref{th:type-of-syzygies}. 

\vspace{2mm}
\noindent
(3) As we have shown in the previous section, 
$B=(B_2',B_1)$ determines a homomorphism $A_B : 4E(3) \rightarrow 2E(2) \oplus 2E(1)$ 
uniquely up to automorphisms of $4E(3)$.  Then the pair $(A_B,B)$  defines a complex $M(N_{A_1})$.

\vspace{2mm}
\noindent
(4) Compute the homology $\ker(B)/\im(A_B)$ of $M(N_{A_1})$ in cohomological degree $0$. 
The homology corresponds to the twisted ideal sheaf of a surface in $\P^4$ 
with the expected invariants. 
Smoothness of this surface can be checked by the Jacobian criterion. 

\vspace{2mm}
\noindent
(5) The hyperplane class of the surface has the same type as (\ref{eq:ample-divisor}). 
Indeed, as in Section \ref{sec:adjunction_process}, 
iterating the adjunction process determines the type of the hyperplane class.

\vspace{2mm}
For the fixed $B_1 $, let $A_1$ be an element of  $ \HHom(4\Omega^3(3),2\Omega^2(2))$ with $N_{A_1}=114$ such that 
\begin{itemize}
\item[$(c_1)$] $C_{A_1}$ is smooth and 
\item[$(c_2)$] $C_{A_1}$ does not intersect $X_{B_1}$. 
\end{itemize}
From the reconstruction of $X$ described above, 
we come up with the question: 

\vspace{2mm}
\noindent
``Given $A_1$, can we construct a smooth rational surface by following 
$(1), (2), (3),(4)$ and $(5)$?" 

\vspace{2mm}
\noindent
The answer is yes. Indeed, by computer search,  
we are able to find an $A_1 \in  \HHom(4\Omega^3(3),2\Omega^2(2))$ 
satisfying Conditions $(c_1)$ and $(c_2)$.  
For this $A_1$, there is a smooth rational surface of the same type as $X$.  
So the next question is 

\vspace{2mm}
\noindent
``What happens, if one varies the value of $N_{A_1}$?" 

\vspace{2mm}
\noindent 
In this section,  using the construction described above, we establish 
the existence of the family of smooth rational surfaces 
in $\P^4$ over a finite field for each $114 \leq N_{A_1} \leq 117$ 
and show that there exists a lift to characteristic zero for each family.  
We also establish the existence of a smooth rational surface  for $N_{A_1}=113$. 
In this case, the intersection of $C_{A_1}$ and $X_{B_1}$ is 
no longer empty.  

\begin{remark}
Let $A_1 \in \HHom(4E(3),2E(2))$ and let 
$S_{A_1}$ be the solution space to the homogeneous system 
of linear equations associated with $A_1$. 
Let $U_{B_1}$ be the vector space formed by 
the pairs of two $E$-linear combinations of  columns of $B_1$. 
Then $U_{B_1}$ can be regarded as a subspace of $\HHom(2E(2),3E)$. 
Since $S_{A_1}$ contains $U_{B_1}$ as 
a  twenty-dimensional subspace, 
both the columns of $B_2'$ can be written 
as $E$-linear combinations of columns of $B_1$ if and only if $N_{A_1}=120$, 
so we may exclude this case. 
\end{remark}


Let $V_{A_1}$ be the column space of $A_1$, and let $V_{B_1}$ be 
the row space of $B_1$. 
Then $\P(V_{A_1})$ and $\P(V_{B_1})$ can be both 
embedded as subvarieties into $G=G(2,V)$.  
Let us denote by $Z_{A_1}$ and $Z_{B_1}$ the images of 
$\P(V_{A_1})$ and $\P(V_{B_1})$ respectively.

\begin{lemma}[Ranestad]
Let $A_1$ be an element of $\HHom(4E(3),2E(2))$ satisfying 
$(c_1)$ and $(c_2)$, and 
let $N_{A_1}$ be the minimal number of equations in the homogeneous 
linear system associated with $A_1$. 
If $Z_{A_1}$ and $Z_{B_1}$ intersect in $k$ points,  
then  $N_{A_1} $ is at most $120-k$. 
\label{th:intersection-number}
\end{lemma}
\begin{proof}
 See Lemma~\ref{th:kristian's-lemma1} in Appendix for the proof. 
\end{proof}
\begin{remark}
(i) Let $\mathfrak{F}$ be the family of rational normal curves in $\P(V)$, 
and let $\mathfrak{F}(N_{A_1})$ be the subfamily of $\mathfrak{F}$ formed by 
rational normal curves $C_{A_1}$ satisfying $(c_2)$ such that  
the minimal number of equations in the homogeneous linear system associated 
with $A_1$ equals $N_{A_1}$. Then
$\codim \left(\mathfrak{F}, \mathfrak{F}(N_{A_1})\right)\leq120-N_{A_1}$ 
by Theorem~\ref{th:intersection-number}.  

\vspace{1mm}
\noindent 
(ii) Each vector of $V_{A_1}$ spans a two-dimensional subspace of $V$. 
Similarly, each vector of $V_{B_1}$ spans a two-dimensional subspace. 
Then $v \in V_{A_1}$ and $w \in V_{B_1}$ span the same subspace $U$ of $V$ if and only if 
$Z_{A_1}$ and $Z_{B_1}$ intersect at the point corresponding to $U$. 
\label{th:codimension}
\end{remark}
The following lemma gives a lower bound for $N(A_1)$: 
\begin{lemma}[Ranestad]
Let $A_1$ be an element of $\HHom(4E(3),2E(2))$ satisfying 
$(c_1)$ and $(c_2)$, and let $\Gamma=Z_{A_1} \cap Z_{B_1}$. 
Suppose that $\Gamma$ is finite. 
Then $\Gamma$ consists of at most six points. 
\label{th:lower-bound}
\end{lemma}
\begin{proof}
See Lemma~\ref{th:kristian's-lemma2} in Appendix for the proof. 
\end{proof}

\vspace{2mm}
Let $N_{A_1}=119$ and let $B_2$ be an element of $\HHom(2E(2),3E)$ 
obtained in (2). In examples, $B=(B_1,B_2')$ has the syzygies 
of the following type 
\begin{eqnarray}
\begin{tabular}{c|ccccc}
\hline
0 & 3 & 2 & . & .\\
-1 & . & 2 & 6  & 10  \\
-2 & . & . & 5 & $\ast$
\end{tabular} 
\label{eq:betti-table3}
\end{eqnarray}
Taking a map from $4E(3)$ to $6E(3)$ at random  
and compositing the map to the second map in (\ref{eq:betti-table3}) 
give a homomorphism $A: 4E(3) \rightarrow 2E(2) \oplus 2E(1)$. 
Then $A$ and $B$ define a complex of type (\ref{eq:sequence}). 
The homology $\ker(B)/\im(A)$, however, has rank three. 

Similarly, the map $B$ has syzygies 
of the following type in the case of $N_{A_1}=118$: 
\begin{eqnarray}
\begin{tabular}{c|ccccc}
\hline
0 & 3 & 2 & . & .\\
-1 & . & 2 & 5  & 3  \\
-2 & . & . & 3 & $\ast$
\end{tabular} 
\label{eq:betti-table4}
\end{eqnarray}
The homology of the associated complex in cohomological degree 0 
is not the ideal sheaf of a surface with the desired invariants.  

\begin{proposition}
There exist smooth rational surfaces in $\P^4$ over $\bar{\F}_5$ with $d=12$ and $\pi=13$ 
embedded by one of the following linear systems:
\begin{itemize}
\item[$\mathrm{(i)}$] $\left|12L-\sum_{i_1=1}^2 4E_{i_1}-\sum_{i_2=3}^{11} 3E_{i_2} 
-\sum_{i_3=12}^{14} 2E_{i_3}-\sum_{i_4=15}^{21} E_{i_4}\right|$, 
\item[$\mathrm{(ii)}$] $\left|12L-\sum_{i_1=1}^3 4E_{i_1}-\sum_{i_2=4}^{9} 3E_{i_2} 
-\sum_{i_3=10}^{15} 2E_{i_3}-\sum_{i_4=16}^{21} E_{i_4}\right|$,
\item[$\mathrm{(iii)}$] $\left|12L-\sum_{i_1=1}^4 4E_{i_1}-\sum_{i_2=5}^{7} 3E_{i_2} 
-\sum_{i_3=8}^{16} 2E_{i_3}-\sum_{i_4=17}^{21} E_{i_4}\right|$,
\item[$\mathrm{(iv)}$] $\left|12L-\sum_{i_1=1}^5 4E_{i_1}-\sum_{i_2=6}^{17} 2E_{i_3} 
-\sum_{i_4=18}^{21} E_{i_4}\right|$. 
\end{itemize}
\label{th:surfaces-in-positive-char}
\end{proposition}
\begin{proof}
By random search over $\F_5$, 
we can find an $A_1 \in \HHom(4E(3),2E(2))$ satisfying $(c_1)$ and $(c_2)$ 
for each $114 \leq N_{A_1} \leq 117$: 
\begin{itemize}
\item[(i)] 
{\footnotesize
$\left(\begin{array}{cccc} 
e_0+2e_1+2e_2  &  -e_3  & -2e_1+e_2+2e_3 & e_0+2e_1+2e_2-e_3+2e_4 \\
e_0+2e_1-e_2+e_3  & -e_4 & e_1-2e_2  &     e_3+2e_4 
\end{array}
\right);
$}
\item[(ii)]  
{\footnotesize
$ \left(
 \begin{array}{cccc}
 2e_3 & -e_0+2e_1-2e_2+2e_3-e_4 &  2e_1+2e_2+2e_4  &  e_0+2e_2+e_3-2e_4 \\
 2e_4  & e_0-2e_1-2e_2-2e_3-2e_4 & e_1+2e_2+e_3-e_4  & e_2+2e_3-2e_4    
 \end{array}
 \right);
$
} 
\item[(iii)]  
{\footnotesize
$\left(
\begin{array}{cccc}
 -e_2-2e_3  &   2e_0+e_2-e_3-2e_4 & -e_1-2e_2-e_3+2e_4 & e_0-e_2+e_3-e_4 \\
 2e_2+2e_3-e_4 &  e_0+2e_2-2e_3-e_4 &  e_1+e_2+e_3-2e_4  & e_2+2e_3+e_4    
\end{array}
\right);
$ 
}
\item[(iv)]  
{\footnotesize 
$\left(
\begin{array}{cccc}
 -e_1+e_4  & -2e_0+e_1-e_2-e_4  & e_1+e_2-2e_3+2e_4 & e_0+2e_2-2e_3+e_4\\
-2e_1-2e_2-2e_3+2e_4  & e_0+2e_1-2e_2-e_3  & e_1+e_2+2e_3-2e_4 & -2e_2-e_4        
 \end{array}
 \right). 
$
}
\end{itemize}
A $\Mac$ script for finding these $A_1$'s can be obtained from~\cite{abo}. 
For each matrix, we can show by $(2), (3)$ and $(4)$ 
that there is a smooth rational surface in $\P^4$ with the desired invariants. 
The type of a linear system embedding the surface into $\P^4$ 
can be determined as in Section~\ref{sec:adjunction_process}. 
\end{proof}
\begin{remark}
(i) To find the matrices given in the proof 
of Proposition~\ref{th:surfaces-in-positive-char} more effectively, 
we take the following extra steps:
For two fixed vectors $v_1$ and $v_2$ that are contained in $V_{B_1}$,  
choose a $2 \times 2$ matrix $A_1'$ with linear entries randomly 
to make the augmented matrix $A_1=(v_1,v_2,A_1')$. 
Then  compute  $N_{A_1}$.  In this case, 
$N_{A_1}$ should be less than or equal to 118 by Remark~\ref{th:codimension} (ii). 
If $N_{A_1} \leq 117$, then  take Steps from (1) to (5). 
We repeat this process until a smooth surface is found. 

\vspace{1mm}
\noindent 
(ii) Let $A_1  \in \HHom(4E(3),2E(2))$ with $114 \leq N_{A_1} \leq 117$ 
satisfying $(c_1)$ and $(c_2)$.  
By Proposition \ref{th:surfaces-in-positive-char}, 
we may assume that there are elements of $\HHom(2E(2),3E)$ 
obtained by  random choices of values for $140-N_{A_1}$ parameters 
that give rise to smooth surfaces in $\P^4$ with $d=12$ and $\pi=13$. 
Let $B_2' $ and $B_2'' $ be such elements of $\HHom(2E(2),3E)$. 
Then the corresponding monads are isomorphic if and only if 
$B_2'$ and $B_2''$ differ only by a constant (modulo $U_{B_1}$).  
This is equivalent to the random choices for $B_2'$ and $B_2''$ are the same up to constant.  
It turns out, therefore, that 
the family of smooth surfaces obtained in this way has dimension 
$
(140-N_{A_1}-1)-20=119-N_{A_1}.
$  

\vspace{1mm}
\noindent 
(iii) For each $A_1$ given in the proof of Proposition~\ref{th:surfaces-in-positive-char}, 
we can check that $Z_{A_1}$ and $Z_{B_1}$ intersect in $120-N_{A_1}$ points. 
So, for a general choice, the equality in Lemma~\ref{th:intersection-number} holds, 
and hence the codimension of $\mathfrak{F}(N_{A_1})$ in $\mathfrak{F}$ is expected to be 
$120-N_{A_1}$. 
\label{th:dimension-of-family}
\end{remark}
The proof of Lemma~\ref{th:lower-bound} suggests us the 
existence of a $2 \times 4$ matrix $A_1$ with entries from $V$ such that 
$C_{A_1}$ is smooth and $Z_{A_1}$ intersects $Z_{B_1}$ in more than 
six points  if we allow $C_{A_1}$ to intersect $X_{B_1}$. 

\begin{proposition}
There exists a smooth rational surface in $\P^4$ 
over $\overline{\F}_3$ with $d=12$ and $\pi=13$ embedded by
\[
\left|12L-4E_1-\sum_{i_2=2}^{13} 3E_{i_2} 
-\sum_{i_4=13}^{21} E_{i_4}\right|. 
\]
\label{th:last-family}
\end{proposition}
\begin{proof}
We can find an $A_1 \in \HHom(4E(3),3E(2))$ 
such that $C_{A_1}$ is smooth and $Z_{A_1}$ intersects $Z_{B_1}$ in seven points 
over $\F_3$ by random search: 
\[
 A_1= 
\left(
\begin{array}{cccc}
 -e_4    &     -e_2-e_3+e_4  &      -e_1  & e_0-e_1-e_2+e_3+e_4 \\
 -e_2-e_3+e_4 &  e_0+e_1+e_2+e_3-e_4 & e_2  &  -e_1-e_2+e_3-e_4    
\end{array}
\right). 
\]
In this example, $N_{A_1}=113$.
So the codimension of $\mathfrak{F}(N_{A_1})$ in $\mathfrak{F}$ 
is expected to be $7$.  
For this $A_1$, there is a smooth rational surface $X$ in $\P^4$ 
with the desired invariants. 
The type of a linear system embedding $X$ into $\P^4$ 
can be determined by the adjunction theory (see Section~\ref{sec:adjunction_process}). 
\end{proof}  
\begin{lemma}
Let $A_1$ be a point of $\mathfrak{F}(N_{A_1})$, 
where $\mathfrak{F}(N_{A_1})$ has codimension $120-N_{A_1}$. 
Then there exists a number field $\LL$ and a prime $\mathfrak{p}$ in $\LL$ 
such that the residue field $\cO_{\LL,\mathfrak{p}}/\mathfrak{p}\cO_{\LL,\mathfrak{p}}$ is 
in $\F_p$.  Furthermore, if the surface $X/\F_p$ corresponding to $A_1$ is smooth, 
then the surface $X/\LL$ corresponding to the generic point 
$\Spec \LL \subset \Spec \cO_{\LL,\mathfrak{p}}$ is also smooth.   
\label{th:lift2}
\end{lemma}
\begin{proof}
See Proposition \ref{th:lift} for the proof. 
\end{proof}
\begin{theorem}
There are at least five different families of smooth rational surfaces in $\P^4$ over $\C$ 
with $d=12$ and  $\pi=13$.
\label{th:lift-to-zero} 
\end{theorem}
\begin{proof}
By Lemma \ref{th:lift2}, 
it suffices to show that for each $ 113 \leq N_{A_1} \leq 117$, the subfamily $\mathfrak{F}(N_{A_1})$ 
of $\mathfrak{F}$ has the desired codimension.  
Let $P$ be the set of isomorphic classes of monads of type 
\[
 4 \Omega^3(3) \rightarrow 2\Omega^2(2) \oplus 2\Omega^1(1) 
 \rightarrow 3\cO. 
\]
Let $T_{P,M(A_1)}$ be the Zariski tangent space 
of $P$ at the point corresponding to 
$M(A_1)$ for each $A_1$ given in 
the proof of Propositions \ref{th:surfaces-in-positive-char} 
and \ref{th:last-family}. 
As in the proof of Theorem \ref{th:existence-of-lift},  
we can show that $\dim\left(T_{P,M(A_1)}\right)=38$ for each $A_1$. 
From Remark \ref{th:dimension-of-family} (ii), the dimension of $\mathfrak{F}(N_{A_1})$ 
is therefore 
\[
 \dim(\mathfrak{F}(N_{A_1}))=\dim(T_{P,M(A_1)})-\left(18+(119-N_{A_1})\right)=N_{A_1}-99.  
\] 
Recall that the parameter space of rational normal curves in $\P(V)$ 
has dimension $21$. 
Thus 
\[
 \codim(\mathfrak{F}(N_{A_1}),\mathfrak{F})=21-(N_{A_1}-99)=120-N_{A_1}, 
\]
which completes the proof. 
\end{proof}


\bigskip 
\noindent 
\begin{tabular}{l}
Hirotachi Abo \\
Deptartment of Mathematics \\
Colorado State University \\
Fort Collins, CO, 80523-1874\\
USA
\end{tabular}

\newpage 
\begin{center}
{\Large Appendix to ``Construction of Rational Surfaces of Degree 12 \\
 in Projective Fourspace" }   \\

\vspace{5mm}
{By Kristian Ranestad} 
\end{center}

Throughout this appendix, the ground field $k$ 
is algebraically closed of any characteristic. 

Fix a vector space  $V\cong k^5$, a $3\times 2$ matrix $B_1$ and a $2\times 4$ matrix $A_1$ with entries from $V$.
Let $S_B$ be the subvariety of $\check\P^4=\P(V^*)$ where $B_1$ has rank $1$, 
and let $C_A$ be the subvariety of $\check\P^4$ where $A_1$ has rank $1$. 
Throughout this note, we make the generality assumption on $A_{1}$ and $B_{1}$, that $C_{A}$ and $S_{B}$ do not intersect and they are both smooth, i.e. 
$C_{A}$ is a rational normal quartic curve, while $S_{B}$ is a rational cubic surface scroll.

Consider furthermore a $3\times 2$ matrix $B_2$ and a $2\times 4$ matrix $A_2$ 
with entries from $\bigwedge ^2 V$ and the $3\times 4$ matrix
$$P_{A,B}=B_{1}\circ A_{2}+B_{2}\circ A_{1}$$ with exterior multiplication, so that the entries of $P_{A,B}$ are in $\bigwedge^3V$.
The matrix equation $P_{A,B}=0$ defines linear equations on the coefficients of $A_{2}$ and $B_{2}$.  
We define $N_{A,B}$ to be the rank of this system of linear equations.
To analyse this rank we define two further varieties associated to $A_{1}$ and $B_{1}$. 
By the generality assumption on these matrices every column of $A_1$ and every row of $B_1$ 
has rank $2$, i.e. defines an element in $G=G(2,V)$ or equivalently a 
plane in $\PP(V^*)$. The columns in $A_{1}$ define $3$-secant planes 
to $C_{A}$, while the rows in $B_{1}$ define planes that intersect 
$S_{B}$ in a conic section.  If $V_{A}$ and $V_{B}$ are the column space of $A_{1}$ 
and the row space of $B_{1}$ respectively, 
then the corresponding map  of $\P(V_{A})$ and $\P(V_{B})$ into 
$G\subset \P(\bigwedge^2V)$ is the double embedding.  We let $Z_{A}$ and 
$Z_{B}$ be the images under these double embeddings of $\P(V_{A})$ and 
$\P(V_{B})$, respectively, i.e. $Z_{A}$  is a Veronese $3$-fold, while $Z_{B}$ is a Veronese surface.
The purpose of this appendix is to prove 

\begin{lemma}  
If  $Z_{A_1}$ and $Z_{B_1}$ have exactly 
$r\leq 6$ common points, 
then $$N_{A,B}\leq 120 - r.$$ 
\label{th:kristian's-lemma1}
\end{lemma}  

\begin{proof}
We first show a lemma that explains the restriction on $r$ in the Lemma 
\ref{th:kristian's-lemma1}.

\begin{lemma} Assume that $A_{1}$ and $B_{1}$ satisfies the above generality assumption.    
 If $\Gamma=Z_{A}\cap Z_{B}$ is finite, then it consists of at most $6$ distinct points.
\label{th:kristian's-lemma2}
\end{lemma}

 \begin{proof}  First we prove that $\Gamma$ either consists of at most $6$ points or 
$4$ points in $\Gamma$  lie on a line in both $\P(V_{A})$ and $\P(V_{B})$.

Notice first that $Z_{A}$ itself spans $\P^9$ while
$Z_{B}$ spans $\P^5$.  This $\P^5$ intersects $G$ in the union of $Z_B$
and a plane $P$.  Thus $\Gamma$ thought of as a subscheme of  $\P^3=\P(V_{A})$ is
contained in the four quadrics defined by
restricting the linear forms that vanish on $Z_{B}$ to $Z_{A}$.  Assume 
now that $\Gamma$ consists of at least $7$ distinct points.
If five of them are in a plane, then the conic through these five is a 
fixed curve in all quadrics through $\Gamma$, which means that the 
intersection of $\P^5$ with $Z_{A}$
contains a curve $C$.  But $\P^5$ intersects $G$ in the union of $Z_B$
and a plane, so the curve $C$ must be contained in this plane, i.e. it
must be a conic, the image of a line in $\P^3$.  The intersection of the
plane and $Z_B$ is also a conic, so $\Gamma$ would contain four collinear
points.

On the other hand, if at most three points are in a plane, there is a
unique  twisted cubic through the six points.  If four points lie in a
plane, this curve degenerates into a conic and a line or three lines.
In either case this possibly reducible twisted cubic lies in three quadrics, and
the six points are defined by four quadrics, a contradiction.

Assume that $Z_B$ and $Z_A$ intersect in finite number of points
and four of them lie on a line both in $\P(V_B)$ and $\P(V_A)$.  The 
image of these two lines are two conics in $\P^5$ that clearly lie in 
the plane $P$. 
On the one hand the four planes in $\P(V^*)$ corresponding to the 
four points each intersects
the rational cubic scroll $S_B$ in a conic
and the rational normal curve $C_A$ in three points.  Let $U$ be the 
union of these four planes.

The conic, which is obtained as the intersection
of the plane $P$ and $Z_B$,
corresponds to the line whose underlying vector space
spanned by two rows of $B_{1}$. The determinant of the $2 \times 2$
matrix generated by these two rows
defines a quadric hypersurface $Q$ containing the union $U$ and $S_B$. 
This quadric may have rank $3$ or $4$.
We will prove that 
$Q$ also contains $C_A$ and that therefore $S_B$ and $C_A$ must intersect.

In case $Q$ has rank $4$, $Q$ is a cone with a vertex. Since 
$C_A$ meets each plane in $U$ in three points, $C_{A}$ is contained 
in $Q$ by the Bezout theorem.  If $Q$ has rank $3$, then any plane defined by a 
linear combination of the two rows contains the same line,
which is the directrix of $S_B$.
Assume that $C_A$ is not contained in $Q$.
Then at least two of three points in each plane
are common for all four planes,
because otherwise the number of intersection number
of $Q$ and $C_A$ is going to be more than $8$.
These points, however lie on the directrix of $S_B$, which is a contradiction.
\end{proof}

We return to the proof of Lemma \ref{th:kristian's-lemma1}.

First note that by the generality assumptions the Lemma 
\ref{th:kristian's-lemma1} applies.  Following that notation we let $Z_{A}$ and 
$Z_{B}$ be the images under these double embeddings of $\P(V_{A_{1}})$ and 
$\P(V_{B_{1}})$ respectively.

Consider now the matrix equation 
$P_{A,B}=B_{1}\circ A_{2}+B_{2}\circ A_{1}=0$.  Since the entries in 
$P_{A,B}$ lie in $\bigwedge^3 V$ which has rank $10$, the rank $N_{A,B}$ 
of the linear system of equations in the coefficients is at most 
$3\times 4\times 10=120$.  These equations are 
parametrized by 
\[
 V_{A_{1}}\otimes V_{B_{1}}\otimes \bigwedge^2V. 
\]
In fact, a row $R_{1}$ in $A_{1}$ and a column $L_{1}$ in $B_{1}$, 
the corresponding row $R_{2}$ in $B_{2}$ and column $L_{2}$ in 
$A_{2}$ and a $2$-vector 
$\omega\in \bigwedge^2V$, define an exterior product  $(R_{1}\cdot 
L_{2}+R_{2}\cdot L_{1})\wedge \omega$ which lies in $\bigwedge^5 V\cong k$.  
Assume that the subspaces of $V$ generated by $L_{1}$ and $R_{1}$ coincide, 
then 
\[(R_{1}\cdot 
L_{2}+R_{2}\cdot L_{1})\wedge \left(\bigwedge^2 L_{1}\right)=0
\] 
independent of $L_{2}$ and $R_{2}$. 
Therefore there is one relation among the coefficients of $P_{A,B}$ 
for each point of intersection of 
$Z_{A}$ and $Z_{B}$ in $\P(\bigwedge^2V)$.  Assume now that there are $r\leq 6$ points of intersection 
    $Z_{A}\cap Z_{B}$, and consider their corresponding tensors in  $V_{A_{1}}\otimes V_{B_{1}}\otimes \bigwedge^2V$.  
Notice that these tensors are all pure, so they have natural 
projections on each factor.  In particular, they are linearly 
independent if they are linearly independent in one factor.  In fact 
we end our proof by showing that the $r \leq 6$ points in $Z_{A}\cap 
Z_{B}$ are linearly independent in $\P(\bigwedge^2V)$.

     Since  $Z_{A}$ and $Z_{B}$ are quadratically embedded, no three 
     points on either of them are collinear.  If four points are 
     coplanar, the plane of their span meets both $Z_{A}$ and $Z_{B}$ 
     in a conic.   In the proof of the Lemma \ref{th:kristian's-lemma2} we saw 
     that this is the case only if $C_{A}$ and 
     $S_{B}$ intersect.   If five points in $Z_{A}\cap Z_{B}$ span only a $\P^3$, 
     this $\P^3$ must intersect the Veronese surface $Z_{B}$ in a 
     conic and a residual point, i.e. four of the five points are 
     coplanar as in the previous case.  
     Finally if $Z_{A}\cap Z_{B}$ 
     consists of six points that span a $\P^4$, then this $\P^4$ 
     intersect $Z_{B}$ in rational normal quartic curve, or in two 
     conics.   But the 
     span of $Z_{B}$ intersects $G$ in the union of the plane 
     $P$ and 
     $Z_{B}$, so the $\P^4$ intersects $G$ in a curve of degree 
     $4$ and a line in $P$ that is bisecant to the curve, or the 
     plane $P$ and a conic that meets $P$ in a point.
     Since no four of the six points are coplanar and $Z_{A}\cap Z_{B}$  
     is finite, the intersection 
     of $\P^4$ and $Z_{A}$ cannot contain a curve.
     The six points in $Z_{A}\cap Z_{B}$ considered as points in 
     $\P(V_{A_{1}})$  therefore lie on five quadrics, which have a 
     finite intersection.  If five of 
     them are in a plane, then the conic in the plane through those 
     five points must lie in each quadric, a contradiction.  
     If at most four points are coplanar, there is a 
     unique, possibly reducible twisted cubic curve through the six 
     points.  But then the six points lie on only four quadrics, a 
     contradiction completing our proof.
\end{proof}

\bigskip 
\noindent 
\begin{tabular}{l}
Kristian Ranestad \\
Mathematisk Institutt \\
Universitetet I Oslo  \\
PP.O.Box 1053 \\
N-0316, Oslo \\
Norway
\end{tabular}

\end{document}